\documentstyle[12pt, amscd]{amsart}
\newcommand{\R}{{\Bbb R}}
\newcommand{\g}{{\bold g}}

\newcommand{\newsection}[1]
{\section{#1}\setcounter{theorem}{0}\setcounter{equation}{0}
\par\noindent}

\newtheorem{theorem}{Theorem}

\newtheorem{lemma}[theorem]{Lemma}
\newtheorem{corrollary}[theorem]{Corollary}

\begin{document}

\title[Global Strichartz estimates]
{Global Strichartz Estimates for Nontrapping\\
$\rule{0pt}{20pt}$ Perturbations of the Laplacian}
\maketitle

\begin{center}
$\text{Hart F. Smith}^1 \text{ and Christopher D. Sogge}^2$\\

\vspace{.1in}

$\rule{0pt}{6pt}^1\text{University of Washington, Seattle, WA 98195}$\\

$\rule{0pt}{6pt}^2\,\text{Johns Hopkins University, Baltimore, MD 21218}$
\end{center}


\vspace{.2in}

\newsection{Introduction}
The purpose of this paper is to establish estimates of Strichartz
type, globally in space and time, for solutions to certain
nontrapping, spatially compact perturbations
of the Minkowski wave equation. Precisely, we consider
the following wave equation on the exterior domain $\Omega$
to a compact obstacle, where the spatial dimension is an odd
integer $n\ge 3$.
\begin{equation}\label{1.1}
\begin{cases}
\partial_t^2 u(t,x) - \Delta_{\g} u(t,x) = F(t,x) \,,
\quad (t,x)\in \R\times\Omega\,,
\\
u(0,x)=f(x)\in \dot{H}^{\gamma}_D(\Omega)\,,
\\
\partial_tu(0,x)=g(x)\in \dot{H}^{\gamma-1}_D(\Omega)\,,
\\
u(t,x)=0 \, ,\quad x\in \partial\Omega \,.
\end{cases}
\end{equation}
The operator $\Delta_{\g}$ is assumed to be
the Laplace-Beltrami operator associated to a smooth,
time-independent Riemannian metric $\g(x)$, such
that $\g_{ij}(x)=\delta_{ij}$ for $|x|\ge R\,.$ 
The set $\Omega$ is assumed to be the
complement in $\R^n$ of a smoothly bounded, compact subset
of $|x|<R$, such that $\Omega$
is strictly geodesically concave with respect to $\g$.
The case that $\Omega =\R^n$ is also permitted. Finally,
\nolinebreak the

\vfill\pagebreak

\noindent
geodesic flow on $\Omega$ with respect to $\g$ and normal
reflection on $\partial\Omega$ is assumed
to be nontrapping, in the sense that all geodesics
exit the set $|x|\le R$ within some fixed, finite time.

The proof consists of showing that exponential energy
decay bounds for \eqref{1.1} (for the case of compactly supported data)
allow one to deduce a global Strichartz type estimate for \eqref{1.1} from
knowledge of the same estimate locally in space-time,
together with the global estimate for solutions of the Minkowski
wave equation on $\R^{1+n}\,.$ The local Strichartz estimates for \eqref{1.1}
were established, in the homogeneous case $F=0$, by the authors
in \cite{SS}. By a special case of a 
lemma of Christ and Kiselev \cite{CK},
this allows one to deduce local
Strichartz estimates for the inhomogeneous problem \eqref{1.1}.
We thank T. Tao for pointing out this fact to us, and we 
include the details in this paper for completeness.

We say that $1\le r,s\le 2\le p,q\le\infty$ and
$\gamma$ are admissible if the
following two mixed norm estimates hold.

\noindent
{\bf Local Strichartz estimates.}
{\it
For data $f,g,F$ supported in $|x|\le R$,
for solutions to \eqref{1.1} the following holds,}
\begin{multline}\label{1.2}
\|u\|_{L^p_tL^q_x([0,1]\times\Omega)}+
\sup_{0\le t\le 1}\|u\|_{H^{\gamma}_D(\Omega)}+
\sup_{0\le t\le 1}\|\partial_t u\|_{H^{\gamma-1}_D(\Omega)}\\
\le C\,
\bigl(\,\|f\|_{H^{\gamma}_D(\Omega)}+\|g\|_{H^{\gamma-1}_D(\Omega)}
+\|F\|_{L^r_tL^s_x([0,1]\times\Omega)}\,\bigr)\,.
\end{multline}

\noindent
{\bf Global Minkowski Strichartz estimates.}
{\it
For solutions to \eqref{1.1}, in the case of $\,\Omega=\R^n$ and 
$\bold g_{ij}(x)=\delta_{ij}$, the following holds,
\begin{equation}\label{1.3}
\|u\|_{L^p_tL^q_x(\R^{1+n})}\le C\,
\bigl(\,\|f\|_{\dot{H}^{\gamma}(\R^n)}+\|g\|_{\dot{H}^{\gamma-1}(\R^n)}
+\|F\|_{L^r_tL^s_x(\R^{1+n})}\,\bigr)\,.
\end{equation}
}

\noindent
Additionally, for technical reasons we need assume $p>r$, and
$\gamma\le (n-1)/2$.
For such a set of indices, we show that the same estimate holds
(with different $C\,$) for solutions to the perturbed equation.

\begin{theorem}\label{theorem1.1}
For admissible $p,q,r,s,\gamma,$ the following estimates holds
for solutions to the mixed Cauchy problem \eqref{1.1}
$$
\|u\|_{L^p_tL^q_x(\R\times\Omega)}\le C\,
\bigl(\,\|f\|_{\dot{H}_D^{\gamma}(\Omega)}
+\|g\|_{\dot{H}_D^{\gamma-1}(\Omega)}
+\|F\|_{L^r_tL^s_x(\R\times\Omega)}\,\bigr)\,.
$$
\end{theorem}

The ingredient that allows us to establish Theorem \ref{theorem1.1}
from \eqref{1.3} and \eqref{1.2} is the following decay estimate
for solutions to the homogeneous problem with compact data.

\noindent{\bf Exponential energy decay.}
{\it 
For data $f,g$ supported in $|x|\le R$, and $\beta(x)$ smooth, supported
in $|x|\le R\,,$
there exist $C<\infty$ and $\alpha>0$ such that
for solutions to \eqref{1.1} where $F=0$ the following holds,}
\begin{equation}\label{1.4}
\|\beta u(t,\cdot)\|_{H_D^{\gamma}(\R^n)}
+\|\beta\partial_t u(t,\cdot)\|_{H_D^{\gamma-1}(\Omega)}
\le C\,e^{-\alpha |t|}
\Bigl(\,\|f\|_{H_D^{\gamma}(\Omega)}+\|g\|_{H_D^{\gamma-1}(\Omega)}\,\Bigr)\,.
\end{equation}

The decay estimate \eqref{1.4}, which depends on
the nontrapping assumption and requires that $n\ge 3$ be odd, 
was first established in the
obstacle framework by Taylor \cite{T}.
For a more general discussion of energy decay estimates, we refer to
Lax-Philips \cite{LP}, and Vainberg \cite{V}.
It was pointed out to us by M. Zworski that quasimode constructions
show that global Strichartz estimates cannot hold if the metric
has an elliptic closed geodesic. See, for example, Colin de Verdi\`ere
\cite{CV} and Ralston \cite{R}.

There is a long history to establishing the global estimates \eqref{1.3} for
the Minkowski case, beginning with the original work by Strichartz
in \cite{Str1,Str2} for the conformal case $p=q=\frac{2n+2}{n-1}$,
$r=s=\frac{2n+2}{n+3}$, and $\gamma=\frac 12\,.$
We mention here the subsequent 
work of Genibre-Velo \cite{GV}, Pecher \cite{P}, Kapitanski \cite{K},
Lindblad-Sogge \cite{LS}, Mockenhaupt-Seeger-Sogge \cite{MSS},
and Keel-Tao \cite{KT}.

Certain local Strichartz estimates \eqref{1.2} for the obstacle problem
were established by the authors in \cite{SS}, for the case
$F= 0\,.$ It will be shown in Theorem \ref{theorem3.2} that these imply the
appropriate estimates in the case $F\ne 0$, based on an argument of
Christ and Kiselev \cite{CK}. 
M. Beals \cite{B} has obtained fixed time $L^p-L^q$ estimates for the
obstacle problem, globally in $t$, 
in the case $\g_{ij}=\delta_{ij}$, and for data
vanishing near $\partial\Omega.$ 
It is not known whether the Strichartz estimates hold, even locally,
if $\partial\Omega$ has a point of convexity,
but related eigenfunction estimates
are known to fail in that case,
as observed by Grieser \cite{G}, and generalized by the authors
in \cite{SS}. It also is not known whether \eqref{1.3} holds 
for every set of indices for which \eqref{1.2} holds,
and we do not address that question in this paper.
We do mention that if the following conditions
are met, then the indices are admissible in the sense of this paper,
provided that $q,s'<\frac{2(n-1)}{n-3}$,
$$
\textstyle
\frac 1p+\frac nq=\frac n2 -\gamma =\frac 1r+\frac ns-2\,,
\quad
\frac 1p = \left(\frac{n-1}2\right)\Bigl(\frac 12-\frac 1q\Bigr)\,,
\quad
\frac 1{r'} = \left(\frac{n-1}2\right)\Bigl(\frac 12-\frac 1{s'}\Bigr)\,.
$$

We make some remarks concerning the 
homogeneous Sobolev spaces $\dot{H}^\gamma_D(\Omega)\,.$
We always deal with $\gamma<\frac n2\,,$ so that smooth cutoffs
of functions in $\dot{H}^\gamma(\R^n)$ are contained in the
inhomogeneous Sobolev space $H^\gamma(\R^n)$, and the two norms
are comparable on functions supported in a fixed ball.
For simplicity, if $R$ is large enough so that 
$g_{ij}(x)=\delta_{ij}$ when $|x|>R$, and 
$\partial\Omega\subset\{|x|<R\}\,,$ then we fix
$\beta\in C^\infty_c$ with $\beta(x)=1$ for $|x|\le R$, and define
$$
\|f\|_{\dot{H}^\gamma_D(\Omega)}=
\|\beta f\|_{H^\gamma_D(\tilde\Omega)}+
\|(1-\beta)f\|_{\dot{H}^\gamma(\R^n)}\,,
$$
where $\tilde\Omega$ is a compact manifold with boundary containing 
$\Omega\cap\{\,|x|\le R\,\}\,.$ For $0\le\gamma\le 2$, 
the space $H^\gamma_D(\tilde\Omega)$ is the usual Dirichlet space
satisfying $f|_{\partial\tilde{\Omega}}=0\,$ (when this
makes sense.) For larger $\gamma$,
the additional compatibility conditions
$$
\Delta_{\g}^j f\in H^{\gamma-2j}_D(\tilde\Omega)\,,\qquad 2j\le\gamma\,,
$$
must be satisfied to insure that solutions to the Cauchy problem
remain in $H^\gamma_D(\tilde\Omega)$. The spaces with $\gamma<0$
are defined by duality; for the Strichartz estimates \eqref{1.3},
always $\gamma-1\ge -1\,,$ so that the additional
compatibility conditions are irrelevant for negative indices.

\newsection{Global Strichartz Estimates}
In this section we provide the proof of Theorem \ref{theorem1.1}, based on
the assumptions \eqref{1.2}, \eqref{1.3}, and \eqref{1.4}. 
By dilating, we will take $R=\frac 12.$

\begin{lemma}\label{lemma2.1}
Let $u$ solve the Cauchy problem \eqref{1.1} with $F$ replaced by $F+G$,
where the data $f,g$ are supported in $\{\,|x|\le 1\,\}$, and
$F, G$ are supported in $\{\,|t|\le 1\,\}\times\{\,|x|\le 1\,\}\,.$
Then, for any $\rho<\alpha$, and any admissible $p,q,r,s,\gamma,$
there exists $C<\infty$ such that the following holds.
\begin{multline*}
\left\|\,e^{\rho(|t|-|x|)}\,u\right\|_{L^p_tL^q_x(\R\times\Omega)}\\
\le C\left(
\|f\|_{H_D^{\gamma}(\Omega)}+\|g\|_{H_D^{\gamma-1}(\Omega)}
+\|F\|_{L^r_tL^s_x(\R\times\Omega)}
+\int\|G(t,\cdot\,)\|_{H_D^{\gamma-1}(\Omega)}\,dt
\right)\,.
\end{multline*}
\end{lemma}
\noindent{\it Proof.}
We establish the estimate on $t\ge 0$. 
Observe that, by (1.2) and Duhamel's principle,
the inequality holds for the $L^p_tL^q_x$ norm of $u$ over 
$[0,1]\times\Omega\,.$ Also, by (1.2), 
\begin{multline*}
\|u(1,\cdot\,)\|_{H_D^{\gamma}(\Omega)}+
\|\partial_t u(1,\cdot\,)\|_{H_D^{\gamma-1}(\Omega)}
\\
\le C\,\left(
\,\|f\|_{H_D^{\gamma}(\R^n)}+\|g\|_{H_D^{\gamma-1}(\R^n)}
+\|F\|_{L^r_tL^s_x(\R\times\Omega)}
+\int\|G(t,\cdot\,)\|_{H_D^{\gamma-1}(\Omega)}
\,\right)\,.
\end{multline*}
By considering $t\ge 1$, we may thus take $F=G=0\,,$
with $f,g$ now supported in $\{\,|x|\le 2\,\}\,.$

We decompose $u=\beta u+(1-\beta)u\,,$ where $\beta(x)=1$ for
$|x|\le \frac 12\,,$ and $\beta(x)=0$ for $|x|\ge 1\,.$
First consider $\beta u\,.$ We write
$$
\partial_t^2(\beta u)-\Delta_{\g}(\beta u)=
\sum_{j=1}^n b_j(x)\partial_{x_j} u + c(x)u\equiv \tilde G(t,x)\,,
$$
where $b_j$ and $c$ are supported in $\frac 12\le|x|\le 1\,,$
and in particular vanish near $\partial\Omega\,.$
By \eqref{1.4} the following holds,
\begin{multline}\label{2.1}
\|\tilde G(t,\cdot\,)\|_{H_D^{\gamma-1}(\Omega)}+
\|\beta u(t,\cdot\,)\|_{H_D^{\gamma}(\Omega)}+
\|\partial_t(\beta u)(t,\cdot\,)\|_{H_D^{\gamma-1}(\Omega)}\\
\le C\,e^{-\alpha t}\,
\Bigl(
\,\|f\|_{H_D^{\gamma}(\Omega)}+\|g\|_{H_D^{\gamma-1}(\Omega)}\Bigr)\,.
\end{multline}
By \eqref{1.2} and Duhamel's principle, it follows that
$$
\|\beta u\|_{L^p_tL^q_x([j,j+1]\times\Omega)}\le C\,e^{-\alpha j}\,
\Bigl(
\,\|f\|_{H_D^{\gamma}(\Omega)}+\|g\|_{H_D^{\gamma-1}(\Omega)}\Bigr)\,,
$$
which easily yields the desired estimate for $\beta u\,.$

Next consider $(1-\beta)\,u\,.$ Since $\Delta_{\g}=\Delta$
on the support of $(1-\beta)u$, we have
$$
\partial_t^2u-\Delta u=-\tilde G\,,
$$
and by Duhamel's principle we have
$$
u(t,x)=u_0(t,x)+\int_0^\infty u_s(t,x)\,ds\,,
$$
where $u_0$ is the solution of the Minkowski wave equation
on $\R^{1+n}$ with
initial data $\bigl(\,(1-\beta)f,(1-\beta)g\,\bigr)$,
and where $u_s(t,x)$ is the
solution of the Minkowski
wave equation on the set $t>s$ with Cauchy data
$\bigl(0,\tilde G(s,\cdot\,)\bigr)$ on the surface $t=s\,.$
(Recall that $\tilde G$ and $(1-\beta)$ vanish
near $\partial\Omega\,.)$ By Huygen's
principle, on the support of $u_s(t,x)$ we have $t\ge s$ and
$t-|x|\in [s-1,s+1]\,,$
so that by \eqref{1.3} and \eqref{2.1} we have
$$
\left\|\,e^{\rho(|t|-|x|)}\,u_s\right\|_{L^p_tL^q_x(\R_+\times\Omega)}\le
C\,e^{(\rho-\alpha)s}\left(
\,\|f\|_{H_D^{\gamma}(\Omega)}+\|g\|_{H_D^{\gamma-1}(\Omega)}
\,\right)\,,
$$
which leads to the desired estimate for $u\,.$\qed

\bigskip

\begin{lemma}\label{lemma2.2}
Let $\beta(x)$ be smooth and supported in $\{|x|\le 1\}\,,$ and
$\;2\gamma\le n-1\,.$ Then the following holds
$$
\int_{-\infty}^\infty\left\|\,\beta(\,\cdot\,)
\left(e^{it|D|}f\right)(t,\cdot\,)\right\|^2_{H^\gamma(\R^n)}\,dt
\le C_{n,\gamma,\beta}\,\|f\|^2_{\dot{H}^\gamma(\R^n)}\,.
$$
\end{lemma}
\noindent{\it Proof.}
By Plancherel's theorem over $t,x\,,$ the left hand side can be written as
$$
\int_0^\infty\int\left|\,
\int\widehat\beta(\xi-\eta)\,\widehat f(\eta)\,\delta(\tau-|\eta|\,)\,d\eta\,
\right|^2\,(1+|\xi|^2)^\gamma\,d\xi\,d\tau\,.
$$
We next apply the Schwarz inequality over $\eta$ to bound this by
\begin{multline*}
\int_0^\infty\int
\left[\,
\int\bigl|\widehat\beta(\xi-\eta)\bigr|\,\delta(\tau-|\eta|\,)\,d\eta
\,\right]\;
\left[\int\bigl|\widehat\beta(\xi-\eta)\bigr|
\,\bigl|\widehat f(\eta)\bigr|^2\,\delta(\tau-|\eta|\,)\,d\eta\,
\right]\\
\times(1+|\xi|^2)^\gamma\,d\xi\,d\tau\,.
\end{multline*}
From the fact that
$$
\sup_\xi\, (1+|\xi|^2)^\gamma\,
\left[\,
\int\bigl|\widehat\beta(\xi-\eta)\bigr|\;\delta(\tau-|\eta|\,)\,d\eta
\,\right]\le
C_{n,\gamma,\beta}\,\text{min}\,\bigl[\tau^{n-1}\,,\,(1+\tau^2)^\gamma\bigr]
$$
this is in turn bounded by
$$
C_{n,\gamma,\beta}\,
\int
\bigl|\widehat f(\eta)\bigr|^2\,
\text{min}\,\bigl[\,|\eta|^{n-1}\,,\,(1+|\eta|^2)^\gamma\bigr]\,d\eta\,\le
C_{n,\gamma,\beta}\,
\|f\|^2_{\dot{H}^\gamma(\R^n)}\,.\qed
$$

\bigskip

\begin{corrollary}\label{corrollary2.3}
Let $\beta$ be a smooth function, supported in $\{\,|x|\le 1\,\}\,.$
Suppose that the global Minkowski Strichartz estimate \eqref{1.3} holds,
and that $2\gamma\le n-1\,.$ Let $u$ solve the Cauchy problem for
the Minkowski wave equation, with data $f,g,F.$ Then the following
holds,
$$
\sup_{|\alpha|\le 1}\;\int_{-\infty}^{\infty}
\bigl\|\,\beta\,\partial_{t,x}^\alpha u\,\bigr\|^2_{H^{\gamma-1}(\R^n)}\,dt
\le C\,\bigl(\,\|f\|_{\dot{H}^\gamma(\R^n)}+\|g\|_{\dot{H}^{\gamma-1}(\R^n)}
+\|F\|_{L^r_tL^s_x(\R^{1+n})}\,\bigr)^2\,.
$$
\end{corrollary}
\noindent{\it Proof.}
If $F=0$, this is a direct consequence of Lemma \ref{lemma2.2} above.
If $f=g=0\,,$ then the Strichartz estimate \eqref{1.3}, duality,
and Huygen's principle imply the following (for $t>0$)
$$
\sup_{|\alpha|\le 1}\;
\bigl\|\,\beta\,\partial_{t,x}^\alpha u(t,\cdot\,)
\,\bigr\|^2_{H^{\gamma-1}(\R^n)}
\le C\,\|F\|_{L^r_tL^s_x(\Gamma_t)}^2\,,
$$
where 
$$
\Gamma_t=\{\,(t',x)\,:\,t'\ge 0\,,\,t'+|x|\in [t-1,t+1]\,\}\,.
$$
Since $r,s\le 2\,,$ the following holds
$$
\int_0^\infty \|F\|_{L^r_tL^s_x(\Gamma_t)}^2\le 
2\,\|F\|_{L^r_tL^s_x(\R^{1+n})}^2\,.\qed
$$

\bigskip

\noindent{\it Proof of Theorem \ref{theorem1.1}}.
By Lemma \ref{lemma2.1}, we may without loss of generality assume that
$f$ and $g$ vanish for $|x|\le 1\,.$
We write 
$$
u=u_0-v=(1-\beta)u_0+\beta u_0-v\,,
$$ 
where $u_0$ solves the Cauchy problem for the Minkowski wave
equation, with data $f,g,F\,,$ where we set $F=0$ on
$\R^n\backslash\Omega\,.$
By \eqref{1.3}, we can restrict attention to $\beta u_0-v\,.$ We write
$$
\bigl(\,\partial_t^2-\Delta_{\g}\,\bigr)\,(\beta u_0-v)=
\beta F+G
$$
where $G(t,x)=\sum_{j=1}^n b_j(x)\partial_{x_j} u_0(t,x) + c(x)u_0(t,x)$
vanishes for $|x|\ge 1\,,$ and satisfies
$$
\int_{-\infty}^\infty\|G(t,\cdot\,)\|^2_{H_D^{\gamma-1}(\Omega)}\,dt
\le
C\,
\Bigl(\,\|f\|_{\dot{H}^{\gamma}(\R^n)}+\|g\|_{\dot{H}^{\gamma-1}(\R^n)}
+\|F\|_{L^r_tL^s_x(\R\times\Omega)}\Bigr)^2\,,
$$
by Corrollary \ref{corrollary2.3}.
Note that the initial data of $\beta u_0-v$ vanishes.
Let $F_j, G_j$ denote the
restrictions of $F,G$ to the set $t\in [j,j+1]\,,$
and write (for $t>0$)
$$
\beta u_0-v=\sum_{j=0}^\infty u_j(t,x)\,,
$$
where $u_j(t,x)$ is the forward solution to 
$\partial_t^2u_j-\Delta_{\g} u_j=\beta F_j+G_j\,.$

By Lemma \ref{lemma2.1}, the following holds
$$
\left\|\,e^{\rho(t-j-|x|)}\,u_j\right\|_{L^p_tL^q_x(\R\times\Omega)}\le
C\,\left(
\|\beta F_j\|_{L^r_tL^s_x(\R\times\Omega)}
+\int_j^{j+1}\,\|G(t,\cdot\,)\|_{H_D^{\gamma-1}(\Omega)}\,dt
\,\right)\,.
$$
Furthermore, $u_j(t,x)$ is supported in the region $t-j-|x|\ge -1\,.$
Consequently, we have
$$
\begin{array}{rcl}
\|\,u\,\|^2_{L^p_tL^q_x(\R\times\Omega)}&\le&
\displaystyle
C\,
\sum_{j=0}^\infty\;\bigl\|\,e^{\rho(t-j-|x|)}u_j\,
\bigr\|_{L^p_tL^q_x(\R\times\Omega)}^2
\\
\\
&\le&
\displaystyle
C\sum_{j=0}^\infty 
\|\,F_j\,\|^2_{L^r_tL^s_x(\R\times\Omega)}
+C\sum_{j=0}^\infty
\left(\int_j^{j+1}\,\|\,G(t,\cdot\,)\,\|_{H_D^{\gamma-1}(\Omega)}
\,dt\right)^2\\
\\
&\le&
\displaystyle
C\|\,F\,\|^2_{L^r_tL^s_x(\R\times\Omega)}
+
C\int_0^\infty\,\|\,G(t,\cdot\,)\,\|^2_{H_D^{\gamma-1}(\Omega)}\,dt
\,,
\end{array}
$$
where we use the fact that $p,q\ge 2\ge r,s\,.$\qed

\bigskip

\newsection{Homogeneous estimates imply inhomogeneous estimates}
In this section, we provide an elementary proof that
mixed-norm estimates for the homogeneous Cauchy problem imply
the appropriate estimates for the inhomogeneous Cauchy problem. 
The proof is quite general, but we present it here in the context
of the obstacle problem. The main ingredient is a
special case of a lemma of Christ
and Kiselev \cite{CK}, the proof of which we present
here for completeness.
We thank T. Tao for pointing out the relevance of this lemma
to inhomogeneous estimates for the wave equation.

\begin{lemma}\label{lemma3.1}
Let $X$ and $Y$ be Banach spaces and assume that $K(t,s)$ is
a continuous function taking its values in $B(X,Y)$, 
the space of bounded linear mappings from
$X$ to $Y$. Suppose that $-\infty\le a<b\le\infty$, and set
$$
Tf(t)=\int_{a}^{b}K(t,s) f(s)\, ds\,.
$$
Assume that
\begin{equation}\label{3.1}
\|Tf\|_{L^q([a,b],Y)}\le C\,\|f\|_{L^p([a,b],X)}.
\end{equation}
Set
$$
W\!f(t)=\int_{a}^t K(t,s)f(s)\, ds.
$$
Then, if $1\le p<q\le \infty\,,$ 
$$
\|\,W\!f\,\|_{L^q([a,b],Y)}\le 
\frac{2^{-2(1/p-1/q)}\cdot 2C}{1-2^{-(1/p-1/q)}}\;
\|f\|_{L^p([a,b],X)}.
$$
\end{lemma}

As remarked in \cite{CK}, if $K(s,t)=1/(t-s)$, then Lemma \ref{lemma3.1}
does not hold in the case $p=q\in (1,\infty)$.

Using this lemma, we establish the following.

\begin{theorem}\label{theorem3.2}
Suppose that the following estimates hold for solutions to the
Cauchy problem \eqref{1.1}, where $F=0\,,$ and the Cauchy data $f,g$ are
supported in the set $|x|\le R\,.$
$$
\begin{array}{rcl}
\|u\|_{L^p_tL^q_x([0,1]\times\Omega)}
&\le& C\,
\bigl(\,\|f\|_{H^{\gamma}_D(\Omega)}+\|g\|_{H^{\gamma-1}_D(\Omega)}\bigr)\,,\\
\\
\|u\|_{L^{r'}_tL^{s'}_x([0,1]\times\Omega)}
&\le& C\,
\bigl(\,\|f\|_{H^{1-\gamma}_D(\Omega)}+\|g\|_{H^{-\gamma}_D(\Omega)}\bigr)\,.
\end{array}
$$
Then the following estimate holds for solutions to the 
inhomogeneous Cauchy problem \eqref{1.1},
provided that $f,g$ and $F$ are supported in the set $|x|\le R\,.$
$$
\|u\|_{L^p_tL^q_x([0,1]\times\Omega)}
\le C\,
\bigl(\,\|f\|_{H^{\gamma}_D(\Omega)}+\|g\|_{H^{\gamma-1}_D(\Omega)}+
\|F\|_{L^r_tL^s_x([0,1]\times\Omega)}\bigr)\,.
$$
\end{theorem}
\noindent{\it Proof of Theorem \ref{theorem3.2}.}
By finite propagation velocity, we may replace $\Omega$ by a 
compact manifold $\tilde\Omega$ with boundary.
Let $\Lambda=\sqrt{-\Delta_{\g}}$, so that the spectrum of 
$\Lambda$ is bounded below.
By the assumptions of the theorem and duality, the following hold.
$$
\begin{array}{rcl}
f\rightarrow \Lambda^{-\gamma}\,e^{\pm it\Lambda}\,f & : &
 L^2(\tilde\Omega)\rightarrow L^p_tL^q_x([0,1]\times\tilde\Omega)\,,\\
\\
\displaystyle
F\rightarrow \Lambda^{\gamma-1}\int_0^1\,
e^{\pm it\Lambda}\,F(s,\cdot\,)\,ds
& : & L^r_tL^s_x([0,1]\times\tilde\Omega)\rightarrow L^2(\tilde\Omega)\,.
\end{array}
$$
Let $K(t,s)=\Lambda^{-1}\,\sin\bigl((t-s)\Lambda\bigr)\,.$ An 
application of the addition formula
to $\sin\bigl((t-s)\Lambda\bigr)$ yields the following.
$$
\int_0^1 K(t,s) F(s,\cdot\,)\,ds : 
L^r_tL^s_x([0,1]\times\tilde\Omega)\rightarrow
L^p_tL^q_x([0,1]\times\tilde\Omega)\,.
$$
By Duhamel's principle and Lemma \ref{lemma3.1}, the result follows.
(To satisfy the conditions of Theorem \ref{theorem3.2}, we should properly
consider a smooth truncation of the wave group to finite
frequencies. The estimates are then uniform, and thus
hold in the limit.)\qed

\bigskip

\noindent
{\it Proof of Lemma \ref{lemma3.1}.}
Since the argument for $q=\infty$ is similar to the case
$q<\infty$, we assume for simplicity that $q$ is finite.

We can normalize $f$ so that
$$\|f\|_{L^p([a,b],X)}=1.$$
We may also assume without loss of generality that $f(s)$
is a continuous function (with values in $X$) and that if
$$F(t)=\int_{a}^t\|f(s)\|^p_X\, ds,$$
then $F: {\Bbb R}\to [0,1]$ is a bijection.  Note then
that if $I\subset [0,1]$ is an interval,
then
\begin{equation}\label{3.2}
\|\chi_{F^{-1}(I)}(s)f(s)\|_{L^p({\Bbb R},X)}=|I|^{1/p}.
\end{equation}

Next consider the set of all dyadic subintervals of $[0,1]\,.$
If $I$ and $J$ are two such
subintervals, we say that $I\sim J$ if the following 
hold.  First, $I$ and $J$ must have the
same length, and $I$ must lie to the left of $J$.  We 
also require that $I$ and $J$ be non-adjacent but have 
adjacent parents (i.e., $I\subset I_0$ and 
$J\subset J_0$ where $I_0$ and $J_0$
are adjacent dyadic intervals of twice the length).  

Note then that if $J$ is fixed there are only two 
intervals with $I\sim J$.  Moreover, for almost
every $(x,y)\in [0,1]^2$ with $x<y$ there is a unique pair 
$I$, $J$ with $I\sim J$ and $x\in I$
and $y\in J$.

If we apply this fact using the variables $x=F(s)$ 
and $y=F(t)$ then we conclude that
the following identity holds almost everywhere:
$$
\begin{array}{rcl}
\chi_{\{(s,t)\in [a,b]^2: s<t\}}(s,t) 
&=& \displaystyle\chi_{\{(x,y)\in [0,1]^2: x<y\}}(x,y)
\\
\\
&=& \displaystyle\sum_{\{I,J: I\sim J\}}\chi_I(x)\chi_J(y)
\\
\\
&= & \displaystyle
\sum_{\{I,J: I\sim J\}}\chi_{F^{-1}(I)}(s)\chi_{F^{-1}(J)}(t).
\end{array}
$$
Consequently,
$$
W\!f=\sum_{\{I,J: I\sim J\}}\chi_{F^{-1}(J)}T(\chi_{F^{-1}(I)}f)\,.
$$
From this we conclude that
\begin{equation}\label{3.3}
\|W\!f\|_{L^q([a,b],Y)}\le \sum_{j=2}^\infty\;
\bigl\|\sum_{\{I,J: I\sim J, \, |I|=2^{-j}\}} \chi_{F^{-1}(J)}
T(\chi_{F^{-1}(I)}f)\bigr\|_{L^q({[a,b]},Y)}\,.
\end{equation}
Since for every $J$ there are at most two $I$ with $I\sim J$,
and the $J$ with $|J|=2^{-j}$ are disjoint, we conclude that
\begin{multline*}
\bigl\|\sum_{\{I,J: I\sim J, \, |I|=2^{-j}\}} \chi_{F^{-1}(J)}
T(\chi_{F^{-1}(I)}f)\bigr\|_{L^q([a,b],Y)}\\
\le 2\left(\sum_{\{I: |I|=2^{-j}\}}
\bigl\|\,T(\chi_{F^{-1}(I)}f)\bigr\|^q_{L^q([a,b],Y)}\right)^{1/q}\,.
\end{multline*}
By \eqref{3.1} and \eqref{3.2}, it follows that
this quantity is bounded by
\begin{multline*}
2C\,\left(\sum_{\{I: |I|=2^{-j}\}}
\|\,\chi_{F^{-1}(I)} f\,\|^q_{L^p([a,b],X)}\right)^{1/q}\\
\le 2C\,
\left(\sum_{\{I: |I|=2^{-j}\}}
2^{-jq/p}
\right)^{1/q}=2^{-j(1/p-1/q)}\cdot 2C\,.
\end{multline*}
The lemma follows by \eqref{3.3} after summing over $j$, 
where we use the fact that $p<q$.\qed


\bigskip

\begin{thebibliography}{MM}

\bibitem{B} M. Beals:
{\em Global time decay of the amplitude of a reflected wave},
Progress in Nonlinear Differential Equations and their
Applications, {\bf\underline{21}}, (1996), 25--44.

\bibitem{CK} M. Christ, A. Kiselev:
{\em Maximal Inequality}, preprint.

\bibitem{CV} Y. Colin de Verdi\`ere:
{\em Quasimodes sur les vari\'et\'es Riemanniennes},
Invent. Math., {\bf\underline{43}}, (1977), 15--52.

\bibitem{GV} J. Genibre, G. Velo:
{\em Generalized Strichartz inequalities for the wave equation},
J. Funct. Anal., {\bf\underline{133}}, (1995), 50--68.

\bibitem{G} D. Grieser:
{\em $L^p$ bounds for eigenfunctions and spectral projections of the
Laplacian near concave boundaries},
Thesis, UCLA, 1992.

\bibitem{K} L. Kapitanski:
{\em Some generalizations of the Strichartz-Brenner inequality},
Leningrad Math J., {\bf\underline{1}}, (1990), 693--726.

\bibitem{KT} M. Keel, T. Tao:
{\em Endpoint Strichartz Estimates},
Amer J. Math., {\bf\underline{120}}, (1998), 955--980.

\bibitem{LP} P. D. Lax, R. S. Philips:
{\em Scattering Theory (Revised Edition)},
Academic Press Inc., 1989.

\bibitem{LS} H. Lindblad, C.D. Sogge:
{\em On existence and scattering with minimal regularity
for semilinear wave equations},
J. Funct. Anal, {\bf\underline{130}}, (1995), 357--426.

\bibitem{MSS} G. Mockenhaupt, A. Seeger, C. D. Sogge:
{\em Local smoothing of Fourier integral operators and Carleson-Sj\"olin
estimates}, J. Amer. Math. Soc., {\bf\underline{6}}, (1993), 65--130.

\bibitem{P} H. Pecher:
{\em Nonlinear small data scattering for the wave and
Klein-Gordan equations},
Math Z., {\bf\underline{185}}, (1984), 261--270.

\bibitem{R} J. Ralston:
{\em Approximate eigenfunctions of the Laplacian},
J. Diff. Geom., {\bf\underline{12}}, (1977), 87--100.

\bibitem{SS} H. Smith, C.D. Sogge:
{\em On the critical semilinear wave equation outside convex obstacles},
J. Amer. Math. Soc., {\bf\underline{8}}, (1995), 879--916.

\bibitem{Str1} R. Strichartz:
{\em A priori estimates for the wave equation and some applications},
J. Funct. Analysis, {\bf\underline{5}}, (1970), 218--235.

\bibitem{Str2} R. Strichartz:
{\em Restriction of Fourier transform
to quadratic surfaces and decay of solutions to the wave equation}
Duke Math. J., {\bf\underline{44}}, (1977), 705--714.

\bibitem{T} M. Taylor:
{\em Grazing rays and reflection of singularities of solutions
to wave equations},
Comm. Pure Appl. Math., {\bf\underline{29}}, (1976), 1--38.

\bibitem{V} B. R. Vainberg:
{\em On the short wave asymptotic behavior of solutions
of stationary problems and the asymptotic behavior as
$t\rightarrow\infty$ of solutions of non-stationary problems},
Russian Math Surveys, {\bf\underline{30:2}}, (1975), 1--55.

\end{thebibliography}
\end{document}